\documentclass{article}
\usepackage{amsmath, amsthm, amsfonts, amssymb}
\usepackage{hyperref}
\usepackage{listings}
\usepackage{xcolor}
\usepackage[utf8]{inputenc}
\usepackage[T1]{fontenc}

\title{An Unexpected Connection Between the Discrete Zeta Function and the Erd\H{o}s--Straus Conjecture Under Mballa's Conjecture}
\author{
  Philemon Urbain Mballa \\
  \small \texttt{philemon-urbain.mballa@etu.u-paris.fr} \\
  \small \texttt{philemonmballa@gmail.com}
}
\date{May 11, 2025}

\begin{document}

\maketitle

\begin{abstract}
In this article, we establish an additive decomposition of the discrete zeta function (for \( s \in \mathbb{N}^*,\ s > 1 \)), more precisely of the function \( 4(\zeta(s) - 1) \), as a series whose general term is of the form \( \frac{1}{x_n(s)} + \frac{1}{y_n(s)} + \frac{1}{z_n(s)} \), where \( x_n(s), y_n(s), z_n(s) \) are solutions of the Erd\H{o}s--Straus conjecture under a personal conjecture (which I will refer to here as Mballa's Conjecture) that I formulated by parametrization in the article: \url{http://arxiv.org/abs/2502.20935}. This connection thus builds a bridge between analysis and Egyptian fractions in general and the Erd\H{o}s conjecture in particular.
\end{abstract}

\section{Introduction}

The Erd\H{o}s--Straus conjecture posits that for all integers \( n \geq 2 \), the rational number \( \frac{4}{n} \) can be written as the sum of three unit fractions:
$$
\frac{4}{n} = \frac{1}{x} + \frac{1}{y} + \frac{1}{z}, \quad x, y, z \in \mathbb{N}^*.
$$
This equation is equivalent to
$$
\frac{4x - n}{nx} = \frac{y + z}{yz}.
$$
In order for the left-hand side to belong to \( \mathbb{N}^* \), it is necessary that
$$
x \geq \left\lfloor \frac{n}{4} \right\rfloor + 1,
$$
where \( \left\lfloor x \right\rfloor \) denotes the integer part of \( x \), and this notation will be used throughout the rest of the article.
In this work, we propose a parameterization of the solutions, which allows us to define a discrete zeta-type function whose terms are derived directly from the decomposition
$$
\frac{4}{n^s} = \frac{1}{x_n(s)} + \frac{1}{y_n(s)} + \frac{1}{z_n(s)}
$$
under Mballa's conjecture. This leads to a functional identity that coincides with \( 4(\zeta(s) - 1) \), and opens a new analytic approach to the problem.

\section{Mballa's Conjecture}

\textbf{Mballa’s Conjecture} posits that For every \( n \in \mathbb{N} \) with \( n \geq 2 \), there exist \( (x, t) \in (\mathbb{N}^*)^2 \) such that
$$
t^2(4x - n)^2 - 2nxt = q^2, \quad \text{with } q \in \mathbb{N}.
$$

This conjecture has been numerically verified in a previous article up to a reasonable bound on \( n \) (see \url{http://arxiv.org/abs/2502.20935} for more details). In the appendix of this paper, we provide a code that will allow the reader, if desired, to test Mballa’s Conjecture numerically.

In the paper \url{http://arxiv.org/abs/2502.20935}, we established that for all \( n \in \mathbb{N} \), with \( n \geq 2 \), if there exists \( x \in \mathbb{N}^* \) such that \( x \geq \left\lfloor \frac{n}{4} \right\rfloor + 1 \), and \( t \in \mathbb{N}^* \), satisfying
$$
t^2(4x - n)^2 - 2nxt = q^2, \quad \text{with } q \in \mathbb{N},
$$
then we can construct
$$
y = t(4x - n) - q, \quad z = t(4x - n) + q \in \mathbb{N}^*
$$
such that
$$
\frac{4}{n} = \frac{1}{x} + \frac{1}{y} + \frac{1}{z}.
$$

It follows from the above that Mballa’s Conjecture implies the Erd\H{o}s--Straus Conjecture.

\section{Construction of the Functional Identity}

For each natural number \( n \geq 2 \), let \( s \in \mathbb{N} \) be fixed with \( s > 1 \), so that \( n^s \in \mathbb{N} \). Under Mballa’s Conjecture, there exist \( (x_n(s), y_n(s), z_n(s)) \in (\mathbb{N}^*)^3 \) such that
$$
\frac{4}{n^s} = \frac{1}{x_n(s)} + \frac{1}{y_n(s)} + \frac{1}{z_n(s)},
$$
with
$$
x_n(s) > \left\lfloor \frac{n^s}{4} \right\rfloor + 1, \quad y_n(s) = t(4x_n(s) - n^s) - q_n(s), \quad z_n(s) = t(4x_n(s) - n^s) + q_n(s),
$$
where \( q_n(s) \in \mathbb{N} \) satisfies
$$
q_n(s)^2 = t^2(4x_n(s) - n^s)^2 - 2t x_n(s) n^s.
$$
Let us consider the following series, which we denote here by \( \zeta_M(s) \):
$$
\zeta_M(s) := \sum_{n=2}^\infty \left( \frac{1}{x_n(s)} + \frac{1}{y_n(s)} + \frac{1}{z_n(s)} \right),
$$
where \( x_n(s), y_n(s), z_n(s) \in \mathbb{N}^* \) are constructed from the decomposition of \( \frac{4}{n^s} \) under the conjecture described above.

Starting from the equality above and summing over \( n \), we obtain:
$$
 \sum_{n=2}^\infty \frac{4}{n^s} = 4(\zeta(s) - 1) = \zeta_M(s) := \sum_{n=2}^\infty \left( \frac{1}{x_n(s)} + \frac{1}{y_n(s)} + \frac{1}{z_n(s)} \right),
$$
which defines a functional identity between the discrete series \( \zeta_M(s) \) and the classical discrete Riemann zeta function.

\section{Convergence of \(\zeta_M(s)\)}

\subsection{Method 1: Direct Equivalence}
We can clearly see that the general term
$$
\frac{1}{x_n(s)} + \frac{1}{t(4x_n(s) - n^s) - q_n(s)} + \frac{1}{t(4x_n(s) - n^s) + q_n(s)} = \frac{4}{n^s},
$$
with \( q_n(s) \) satisfying
$$
t^2(4x_n(s) - n^s)^2 - 2tn^s x_n(s) = q_n(s)^2.
$$
Now, the term \( \frac{4}{n^s} \) belongs to a multiplicatively adjusted Riemann-type series, which converges for \( s > 1 \). Therefore, the series \( \zeta_M(s) \) converges.

\subsection{Method 2: Comparison Theorem}

\textbf{Proposition (Linearity of Series).} Let \( (u_n) \) and \( (v_n) \) be two sequences with values in a field \( \mathbb{K} \), and let \( \lambda \in \mathbb{K} \). If the series \( \sum u_n \) and \( \sum v_n \) converge, then the series \( \sum (u_n + \lambda v_n) \) also converges, and we have
$$
\sum (u_n + \lambda v_n) = \sum u_n + \lambda \sum v_n.
$$

Let us apply this proposition to prove the convergence of the series \( \zeta_M(s) \).

\subsection*{Convergence of the Series \(\sum_{n=2}^{\infty} \frac{1}{x_n(s)}\)}

We know that \( x_n(s) \geq \left\lfloor \frac{n^s}{4} \right\rfloor + 1 \), a condition that was analytically imposed in the introduction.
By definition, for any real number \( a > 0 \), we have \( \left\lfloor a \right\rfloor + 1 > a \). Therefore,
\begin{equation} \label{eq:1}
x_n(s) > \frac{n^s}{4} \quad \Rightarrow \quad \frac{1}{x_n(s)} < \frac{4}{n^s}.
\end{equation}

Since the general term \( \frac{4}{n^s} \) defines a Riemann-type series that converges for \( s > 1 \), it follows by the comparison theorem that the series \( \sum_{n=2}^{\infty} \frac{1}{x(s)} \) also converges.

\subsection*{Convergence of the Series \(\sum_{n=2}^{\infty} \left( \frac{1}{y_n(s)} + \frac{1}{z_n(s)} \right)\)}

We observe that:
$$
\left| \frac{1}{y_n(s)} + \frac{1}{z_n(s)} \right| = \left| \frac{1}{t(4x_n(s) - n^s) - q_n(s)} + \frac{1}{t(4x_n(s) - n^s) + q_n(s)} \right|.
$$
By the standard identity of conjugate reciprocals, this yields:
$$
\left| \frac{1}{y_n(s)} + \frac{1}{z_n(s)} \right| = \left| \frac{2t(4x_n(s) - n^s)}{t^2(4x_n(s) - n^s)^2 - q_n(s)^2} \right|.
$$
But since
$$
t^2(4x_n(s) - n^s)^2 - q_n(s)^2 = 2tn^s x_n(s),
$$
we have:
$$
\left| \frac{1}{y_n(s)} + \frac{1}{z_n(s)} \right| = \frac{|4x_n(s) - n^s|}{n^s x_n(s)}.
$$

Now, by the triangle inequality and the positivity of \( x_n(s) \) and \( n^s \), we get:
$$
|4x_n(s) - n^s| \leq 4x_n(s) + n^s \quad \Rightarrow \quad \frac{|4x_n(s) - n^s|}{n^s x_n(s)} \leq \frac{4}{n^s} + \frac{1}{x_n(s)}.
$$
Using inequality~\eqref{eq:1}, we deduce:
$$
\left| \frac{1}{y_n(s)} + \frac{1}{z_n(s)} \right| < \frac{8}{n^s}.
$$

Thus, the series \( \sum_{n=2}^{\infty} \left( \frac{1}{y_n(s)} + \frac{1}{z_n(s)} \right) \) converges (absolutely) by comparison with a convergent Riemann series (for \( s > 1 \)).

\subsection*{Conclusion}

By the linearity of convergent series, we conclude that the series
$$
\zeta_M(s) = \sum_{n=2}^{\infty} \left( \frac{1}{x_n(s)} + \frac{1}{y_n(s)} + \frac{1}{z_n(s)} \right)
$$
also converges.

Therefore, we have the identity:
\begin{equation}
\zeta_M(s) = \sum_{n=2}^\infty \frac{4}{n^s} = 4(\zeta(s) - 1).
\label{eq:2}
\end{equation}

Why is this relation interesting? It expresses the discrete Riemann zeta function (which lies at the heart of mathematics) in terms of the series \( \zeta_M \), under Mballa’s Conjecture, which itself implies the Erd\H{o}s--Straus Conjecture. Therefore, the more information we have about Mballa’s Conjecture (in particular, the structure of \( q_n(s) \), \( x_n(s) \), etc.), the more insight we gain about the Erd\H{o}s--Straus Conjecture. Most importantly, this also provides a new angle for understanding the values of the zeta function at odd integers, which remain one of the major mysteries in mathematics.

Indeed, as stated in \eqref{eq:2}, we have:
$$
4(\zeta(2k + 1) - 1) = \zeta_M(2k + 1), \quad \text{for } k \in \mathbb{N}^*.
$$

\section{Numerical Evidence}

\subsection*{Remark on Numerical Simulations}

Although numerical simulations can be useful to illustrate the accuracy and convergence behavior of the identity
$$
\zeta_M(s) = 4(\zeta(s) - 1),
$$
they are not strictly necessary in our case. This identity follows rigorously from the pointwise equality
$$
\frac{1}{x_n(s)} + \frac{1}{y_n(s)} + \frac{1}{z_n(s)} = \frac{4}{n^s},
$$
which holds for all \( n \geq 2 \) under Mballa’s Conjecture. Given that both series converge, the equality of the infinite sums is ensured by classical results on term-by-term summation.

Nevertheless, for the interested reader, we provide in the appendix a Python script that numerically verifies this identity up to any chosen value of \( n \), for fixed \( s \in \mathbb{N}^* \).

\section*{References}

[1] P.~U.~Mballa, \textit{Partial Resolution of the Erdős–Straus, Sierpiński, and Generalized Erdős–Straus Conjectures Using New Analytical Formulas}, arXiv:2502.20935 [math.NT], 2025. Available at: \url{http://arxiv.org/abs/2502.20935}

[2] P.~Erdős, \textit{On a Diophantine Equation}, Mat. Lapok, vol.~1, pp.~192–210, 1950.

[3] J.~Derbyshire, \textit{The Prime Obsession: Bernhard Riemann and the Greatest Unsolved Problem in Mathematics}, 1st ed., Walker \& Company, 2003.

\appendix
\section*{Appendix: Numerical Verification Code}

The following Python script allows the reader to verify the identity
$$
\zeta_M(s) = 4(\zeta(s) - 1)
$$
numerically up to any chosen value \( N_{\text{max}} \), for fixed \( s \in \mathbb{N}^* \), \( s > 1 \). 

For example, for \( s = 3 \) and \( N_{\text{max}} = 150 \), the program returns:
\begin{quote}
\small

\( \zeta(3) - 1 \approx 0.20203482859170052 \) \\
\texttt{Approach via series: 0.20203482859170052} \\
\texttt{Absolute error: 0.0}
\end{quote}

This shows a perfect numerical match between both sides of the identity, confirming the analytic equality with machine precision.

Below is the full Python script used for the computation:

\lstset{
    language=Python,
    basicstyle=\ttfamily\footnotesize,
    keywordstyle=\color{blue},
    commentstyle=\color{gray},
    stringstyle=\color{red},
    breaklines=true,
    numbers=left,
    numberstyle=\tiny,
    frame=single,
    showstringspaces=false
}

\begin{lstlisting}
import numpy as np

def find_triplet_corrected(n, s=2, a=4):
    """
    Find (x, y, z) in (N^*)^3 such that:
    4 / n^s = 1/x + 1/y + 1/z,
    using:
        y = t(4x - n^s) - q,
        z = t(4x - n^s) + q,
        with q^2 = t^2(4x - n^s)^2 - 2tn^s x
    """
    n_s = n**s
    for x in range(n_s // a + 1, 300 * n_s):
        t_min = max(1, (2 * n_s * x) // ((a * x - n_s)**2))
        for t in range(t_min, t_min + 500):
            q_square = t**2 * (a * x - n_s)**2 - 2 * t * x * n_s
            if q_square >= 0:
                q = int(np.sqrt(q_square))
                if q * q == q_square:
                    y = t * (a * x - n_s) - q
                    z = t * (a * x - n_s) + q
                    if y > 0 and z > 0:
                        return x, y, z
    return None, None, None

def test_series_identity_corrected(N_max, s=2):
    """
    Test: sum_{n=2}^{N_max} (1/x + 1/y + 1/z)/4 ~= sum_{n=2}^{N_max} 1/n^s = zeta(s) - 1
    """
    left_sum = 0  # zeta(s) - 1
    right_sum = 0  # Series using parametrization

    for n in range(2, N_max + 1):
        left_sum += 1 / (n**s)
        x, y, z = find_triplet_corrected(n, s)
        if x and y and z:
            right_sum += (1 / x + 1 / y + 1 / z) / 4

    error = abs(left_sum - right_sum)
    print("N_max =", N_max)
    print("zeta(", s, ") - 1 ~= ", left_sum)
    print("Series approach:", right_sum)
    print("Absolute error:", error)

    return left_sum, right_sum, error

# Example
test_series_identity_corrected(N_max=150, s=3)
\end{lstlisting}

\appendix

\section*{Appendix: Python Code for Verifying Mballa’s Conjecture}

The following Python code tests Mballa’s Conjecture numerically for all integers \( n \in [n_{\min}, n_{\max}] \). The conjecture states that for every \( n \geq 2 \), there exist integers \( x > \left\lfloor \frac{n}{4} \right\rfloor \) and \( t \in \mathbb{N}^* \) such that:
$$
t^2(4x - n)^2 - 2nxt = q^2, \quad \text{with } q \in \mathbb{N}.
$$

The search is conducted over the range \( x \in [\lfloor n/4 \rfloor + 1,\ 300n] \), and for each \( x \), the parameter \( t \) is tested in the interval \( [t_{\min},\ t_{\min} + 500] \), where:
$$
t_{\min} = \max\left(1,\ \left\lfloor \frac{2n x}{(4x - n)^2} \right\rfloor\right).
$$

\lstset{
    language=Python,
    basicstyle=\ttfamily\footnotesize,
    keywordstyle=\color{blue},
    commentstyle=\color{gray},
    stringstyle=\color{red},
    breaklines=true,
    numbers=left,
    numberstyle=\tiny,
    frame=single,
    showstringspaces=false
}

\begin{lstlisting}
def verify_mballa_conjecture(n_min=2, n_max=100, a=4):
    """
    Numerically verifies Mballa's Conjecture for n in [n_min, n_max].

    The conjecture states: for every n >= 2, there exist (x, t) in (N^*)^2 such that
    t^2(4x - n)^2 - 2nxt = q^2, with q in N.
    

    This version reports:
    - each successful match (x, t, q),
    - the percentage of values captured,
    - the list of values not captured within search bounds.
    """

    captured = 0
    failed_n_list = []

    for n in range(n_min, n_max + 1):
        found = False
        n_s = n
        for x in range(n_s // a + 1, 300 * n_s):
            denom = (a * x - n_s)**2
            if denom == 0:
                continue
            t_min = max(1, (2 * n_s * x) // denom)
            for t in range(t_min, t_min + 500):
                q_square = t**2 * (a * x - n_s)**2 - 2 * t * x * n_s
                if q_square >= 0:
                    q = int(q_square**0.5)
                    if q * q == q_square:
                        print(f"n = {n}: OK (x = {x}, t = {t}, q = {q})")
                        captured += 1
                        found = True
                        break
            if found:
                break
        if not found:
            print(f"n = {n}:  No solution found within bounds")
            failed_n_list.append(n)

    total = n_max - n_min + 1
    success_rate = 100 * captured / total

    print("\n========== Summary ==========")
    print(f"Interval tested: n in [{n_min}, {n_max}]")
    print(f"Captured: {captured} out of {total} -> {success_rate:.2f}%")
    if failed_n_list:
        print("Failed to capture the following n values:")
        print(failed_n_list)
    else:
        print("Ok All values of n were successfully captured.")

# Example execution
verify_mballa_conjecture(n_min=2, n_max=150)
\end{lstlisting}

\section*{Appendix (continued): Python Code for the Discrete Zeta Version of Mballa’s Conjecture}

The following Python script numerically tests the discrete zeta formulation of Mballa’s Conjecture, where the classical parameter \( n \) is replaced by \( n^s \), with \( s \in \mathbb{N}^*,\ s > 1 \). The generalized equation to be solved is:
$$
t^2(4x - n^s)^2 - 2tn^s x = q^2, \quad \text{with } q \in \mathbb{N}.
$$
For each integer \( n \in [n_{\min}, n_{\max}] \), this script searches for integers \( x > \left\lfloor \frac{n^s}{4} \right\rfloor \) and \( t \in \mathbb{N}^* \) such that the above identity holds.

The search bounds are:
$$
x \in \left[ \left\lfloor \frac{n^s}{4} \right\rfloor + 1,\ 300n^s \right], \quad
t \in \left[ t_{\min},\ t_{\min} + 500 \right], \quad
t_{\min} = \max\left(1,\ \left\lfloor \frac{2n^s x}{(4x - n^s)^2} \right\rfloor \right).
$$

\begin{lstlisting}[language=Python]
def verify_mballa_conjecture_zeta_version(n_min=2, n_max=100, s=2, a=4):
    """
    Verifies the generalized (discrete zeta) version of Mballa's Conjecture.

    For each n >= 2, this version tests whether there exist (x, t) in (N^*)^2 such that:
        t^2(4x - n^s)^2 - 2tn^s x = q^2, with q in N.

    Reports each successful match and the final success rate.
    """

    captured = 0
    failed_n_list = []

    for n in range(n_min, n_max + 1):
        found = False
        n_s = n**s
        for x in range(n_s // a + 1, 300 * n_s):
            denom = (a * x - n_s)**2
            if denom == 0:
                continue
            t_min = max(1, (2 * n_s * x) // denom)
            for t in range(t_min, t_min + 500):
                q_square = t**2 * (a * x - n_s)**2 - 2 * t * x * n_s
                if q_square >= 0:
                    q = int(q_square**0.5)
                    if q * q == q_square:
                        print(f"n = {n}: OK (x = {x}, t = {t}, q = {q}) for s = {s}")
                        captured += 1
                        found = True
                        break
            if found:
                break
        if not found:
            print(f"n = {n}: No solution found within bounds for s = {s}")
            failed_n_list.append(n)

    total = n_max - n_min + 1
    success_rate = 100 * captured / total

    print("\n========== Zeta Version Summary ==========")
    print(f"Interval tested: n in [{n_min}, {n_max}] with s = {s}")
    print(f"Captured: {captured} out of {total} -> {success_rate:.2f}%")
    if failed_n_list:
        print("Failed to capture the following n values:")
        print(failed_n_list)
    else:
        print("All values were successfully captured for the zeta version.")

# Example usage:
verify_mballa_conjecture_zeta_version(n_min=2, n_max=100, s=3)
\end{lstlisting}

\end{document}